\documentclass[12pt]{article}

\usepackage{amssymb,amsmath,amsthm}
\usepackage{amsfonts}
\usepackage{graphicx}
\usepackage{subcaption}
\usepackage{float}
\usepackage{geometry}
\usepackage{hyperref}
\usepackage{xcolor}
\usepackage{algorithm}
\usepackage{algpseudocode}
\usepackage{booktabs}
\usepackage{cite}

\renewcommand{\mathbf}[1]{\boldsymbol{#1}}

\newtheorem{theorem}{Theorem}[section]

\newtheorem{remark}[theorem]{Remark}

\newcommand{\calL}{\mathcal{L}}

\newcommand{\Eg}{\operatorname{div}_g}
\newcommand{\Dg}{\Delta_g}
\newcommand{\Ng}{\nabla_g}
\newcommand{\vth}{{\boldsymbol{\vartheta}}}

\newcommand{\dVol}{\,\mathrm{dVol}_g}
\newcommand{\PiM}{\Pi_M}

\title{\Large\bfseries Neural Pushforward Samplers for the Fokker--Planck
  Equation on Embedded Riemannian Manifolds}

\author{Andrew Qing He\thanks{Department of Mathematics, Southern Methodist
    University, Dallas, TX 75275.  \texttt{andrewho@smu.edu}}
  \and
  Wei Cai\thanks{Department of Mathematics, Southern Methodist University,
    Dallas, TX 75275.  \texttt{cai@smu.edu}}}

\date{March 2026}

\begin{document}
\maketitle

\begin{abstract}
In this paper, we extend the weak adversarial neural pushforward mapping method to the Fokker--Planck equation on compact embedded Riemannian
manifolds.
The method represents the solution as a probability distribution via a
neural pushforward map that is constrained to the manifold by a retraction
layer, enforcing manifold membership and probability conservation by
construction.
Training is guided by a weak adversarial objective using ambient plane-wave
test functions, whose intrinsic differential operators are derived in closed
form from the geometry of the embedding, yielding a fully mesh-free and
chart-free algorithm.
Both steady-state and time-dependent formulations are developed, and
numerical results on a double-well problem on the two-sphere demonstrate
the capability of the method in capturing multimodal invariant distributions
on curved spaces.
\end{abstract}



\tableofcontents
\newpage

\section{Introduction}
\label{sec:intro}

A wide class of stochastic differential equations (SDEs) arising in molecular
dynamics, robotics, directional statistics, and climate science evolve not in
flat Euclidean space but on a constraint surface: a sphere, a torus, a
Stiefel manifold, or more generally a smooth Riemannian manifold $(M, g)$
\cite{hsu2002stochastic,girolami2011riemann,lei2019geometric}.
The natural companion PDE governing the time evolution of the probability
density is the Fokker--Planck equation (FPE) on $(M, g)$,
\begin{equation}\label{eq:FP_intro}
  \frac{\partial\rho}{\partial t}
  = \frac{\sigma^2}{2}\,\Dg\,\rho - \Eg(\rho\,\mathbf{b}),
\end{equation}
where $\Dg$ is the Laplace--Beltrami operator, $\mathbf{b}$ is a smooth
tangential drift, and $\sigma > 0$ is the diffusion coefficient.

Solving \eqref{eq:FP_intro} numerically is challenging for two independent
reasons.  First, classical finite-difference or finite-element methods require
an explicit triangulation of $M$, which becomes prohibitively expensive in
moderate to high dimension.  Second, even for simple manifolds such as
$S^{n-1}$ with $n \geq 4$, the curse of dimensionality prevents grid-based
methods from being practically viable.

Recent work has addressed these difficulties for flat-space FPEs by
representing the solution \emph{distribution} rather than a pointwise density,
through a neural pushforward map trained with a weak adversarial objective
\cite{He2025WANPF}.  In that approach, a generator network $F_\vth$ pushes a
simple base distribution forward to approximate the FPE solution at any queried
time, with training guided by a three-term Monte Carlo weak formulation that
avoids any spatial mesh.

The present paper extends this Weak Adversarial Neural Pushforward Method
(WANPF) to the Riemannian setting.  The extension rests on two structural
observations.

\begin{enumerate}
  \item \textbf{Non-invertibility is an asset.}
  The pushforward map $F_\vth : \mathbb{R}^d \to M$ need not be invertible,
  and the base dimension $d$ need not equal $\dim M$.
  This is impossible for normalizing-flow or change-of-variables methods, and it
  is what allows WANPF to use a higher-dimensional noise space for greater
  representational capacity.

  \item \textbf{Weak integrals become pointwise expectations.}
  After integrating by parts, every spatial operator in the weak form acts on
  the test function rather than on $\rho$.  The Laplace--Beltrami operator
  applied to an ambient plane-wave test function can be computed in closed
  form at each sample, using only the tangential projection $P(x) = I - xx^\top$
  (for the sphere) and the mean-curvature vector $H(x)$.  No automatic
  differentiation through the intrinsic geometry is required.
\end{enumerate}

The result is a fully mesh-free, chart-free, and Jacobian-free training
procedure.  The generator's samples are constrained to $M$ via a manifold
retraction $\Pi_M$ (e.g., normalization to the unit sphere), and the
initial condition is satisfied automatically by the pushforward architecture.

\paragraph{Contributions.}
\begin{itemize}
  \item A complete steady-state WANPF formulation for the manifold FPE,
        including the adversarial min-max loss and its Monte Carlo estimator
        (Section~\ref{sec:steady}).
  \item A time-dependent WANPF formulation with a three-term weak form
        (Section~\ref{sec:timedep}), matching the flat-space structure
        of \cite{He2025WANPF}.
  \item Closed-form Laplace--Beltrami plane-wave formulae for $S^{n-1}$
        and $\mathbb{T}^n$ (Sections~\ref{sec:pw}--\ref{sec:manifolds}).
  \item Numerical experiments on the double-well FPE on $S^2$, demonstrating
        that the trained generator correctly concentrates mass near the two
        Boltzmann minima (Section~\ref{sec:numerics}).
\end{itemize}

\paragraph{Organization.}
Section~\ref{sec:setting} introduces the manifold FPE and fixes notation.
Section~\ref{sec:steady} derives the steady-state WANPF.
Section~\ref{sec:pw} derives the ambient plane-wave Laplace--Beltrami formulae.
Section~\ref{sec:arch} specifies the pushforward architecture.
Section~\ref{sec:timedep} derives the time-dependent WANPF.
Section~\ref{sec:manifolds} records explicit formulae for the sphere and torus.
Section~\ref{sec:numerics} presents the numerical experiments.
Section~\ref{sec:discussion} discusses extensions.

\section{Setting and Notation}
\label{sec:setting}

Throughout, $(M, g)$ is a smooth, compact, connected Riemannian manifold
without boundary, isometrically embedded in $\mathbb{R}^n$.
Compactness ensures global probability conservation and the absence of spatial
boundary conditions.  Let $m = \dim M \leq n - 1$.

For any $x \in M$, let $T_x M \subset \mathbb{R}^n$ denote the tangent space
of $M$ at $x$.  We write $P(x) \in \mathbb{R}^{n \times n}$ for the orthogonal
projection onto $T_x M$, so $P(x)^2 = P(x)$, $P(x)^\top = P(x)$, and
$P(x) v$ is the tangential component of any ambient vector $v \in \mathbb{R}^n$.

The \emph{mean-curvature vector} $H(x) \in \mathbb{R}^n$ of the embedding is
\begin{equation}\label{eq:H_def}
  H(x) = \operatorname{tr}\!\bigl(\partial_i P(x)\bigr),
\end{equation}
where the trace is over the $n$ coordinate directions.
For a smooth ambient function $f : \mathbb{R}^n \to \mathbb{R}$, the
intrinsic gradient and Laplace--Beltrami operator on $M$ are:
\begin{align}
  \Ng f(x) &= P(x)\,\nabla f(x),
  \label{eq:tang_grad}\\
  \Dg f(x) &= \nabla^2 f(x) : P(x) - \nabla f(x) \cdot H(x),
  \label{eq:LB_ambient}
\end{align}
where $\nabla f$ and $\nabla^2 f$ are the Euclidean gradient and Hessian
of the ambient extension of $f$, and $A : B = \operatorname{tr}(A^\top B)$.
For flat subspaces, $H = 0$ and \eqref{eq:LB_ambient} reduces to the
standard Laplacian restricted to the subspace.

The drift $\mathbf{b}$ is a smooth tangential vector field on $M$, and in
applications derived from a potential $V : \mathbb{R}^n \to \mathbb{R}$ it
takes the form $\mathbf{b}(x) = -P(x)\nabla V(x)$, which is automatically
tangential.

\section{Steady-State WANPF on a Manifold}
\label{sec:steady}

\subsection{Strong Formulation}
\label{sec:steady_strong}

The steady-state FPE seeks a probability measure $\rho$ on $M$ satisfying:
\begin{equation}\label{eq:FP_steady}
  \frac{\sigma^2}{2}\,\Dg\,\rho - \Eg(\rho\,\mathbf{b}) = 0,
  \qquad
  \int_M \rho\dVol = 1, \quad \rho \geq 0.
\end{equation}
On a compact manifold without boundary the Fredholm alternative guarantees
the existence and, under mild irreducibility conditions on $\mathbf{b}$,
the uniqueness of an invariant measure.  When $\mathbf{b} = -\nabla_g V$
for a potential $V$, the unique invariant measure is the Gibbs distribution
$\rho_\infty \propto e^{-2V/\sigma^2}$.

\subsection{Weak Formulation}
\label{sec:steady_weak}

Multiplying \eqref{eq:FP_steady} by a test function $f \in C^\infty(M)$,
integrating over $M$ against $\dVol$, and integrating by parts (boundary
terms vanish by compactness) gives:
\begin{equation}\label{eq:weak_steady}
  \int_M \rho(x)\,\calL_{\mathrm{ss}} f(x)\dVol = 0,
\end{equation}
where the \emph{backward Kolmogorov operator} for the steady-state problem is:
\begin{equation}\label{eq:Lss}
  \calL_{\mathrm{ss}} f
  := \frac{\sigma^2}{2}\,\Dg f + \langle\mathbf{b}, \Ng f\rangle_g.
\end{equation}
Equation \eqref{eq:weak_steady} is an expectation condition:
\begin{equation}\label{eq:steady_expectation}
  \mathbb{E}_{\mathbf{x}\sim\rho}\!\left[\calL_{\mathrm{ss}} f(\mathbf{x})\right] = 0
  \quad \forall \text{ admissible } f.
\end{equation}

\subsection{Pushforward Architecture and Adversarial Loss}
\label{sec:steady_arch}

For the steady-state problem the pushforward map has no temporal argument:
\begin{equation}\label{eq:pf_steady}
  F_{\vth} : \mathbb{R}^d \to M, \qquad
  F_{\vth}(\mathbf{r}) = \PiM\!\bigl(\tilde{F}_{\vth}(\mathbf{r})\bigr),
\end{equation}
where $\tilde{F}_{\vth} : \mathbb{R}^d \to \mathbb{R}^n$ is an unconstrained
MLP and $\PiM$ is the manifold retraction (Section~\ref{sec:arch}).
Given $M_s$ samples $\mathbf{r}^{(m)} \sim \mathcal{P}_{\mathrm{base}}$ and
$K$ adversarial test functions $\{f^{(k)}\}$, the empirical estimator is
\begin{equation}\label{eq:Ehat_steady}
  \hat{E}^{(k)} =
  \frac{1}{M_s}\sum_{m=1}^{M_s} \calL_{\mathrm{ss}} f^{(k)}
  \!\bigl(F_\vth(\mathbf{r}^{(m)})\bigr).
\end{equation}
The adversarial training objective is the min-max problem:
\begin{equation}\label{eq:loss_steady}
  \min_\vth\;\max_{\{\boldsymbol{\eta}^{(k)}\}}\;
  \frac{1}{K}\sum_{k=1}^K \bigl(\hat{E}^{(k)}\bigr)^2,
\end{equation}
where $\boldsymbol{\eta}^{(k)} = \{\mathbf{w}^{(k)}, b^{(k)}\}$ are the
test-function parameters defined in Section~\ref{sec:pw}.

\section{Ambient Plane-Wave Test Functions}
\label{sec:pw}

\subsection{Definition}

Following the flat-space WANPF \cite{He2025WANPF}, the test functions are
ambient plane waves:
\begin{equation}\label{eq:planewave_ss}
  f^{(k)}(x) = \sin\!\bigl(\mathbf{w}^{(k)}\cdot x + b^{(k)}\bigr),
  \qquad k = 1,\ldots,K,
\end{equation}
with learnable parameters $\mathbf{w}^{(k)} \in \mathbb{R}^n$ and
$b^{(k)} \in \mathbb{R}$.  For the time-dependent case an additional
temporal phase $\kappa^{(k)} t$ is included
(Section~\ref{sec:timedep}).

\subsection{Laplace--Beltrami Action on Plane Waves}
\label{sec:pw_LB}

Setting $\phi^{(k)} := \mathbf{w}^{(k)}\cdot x + b^{(k)}$,
the ambient Euclidean derivatives are
$\nabla f^{(k)} = \cos(\phi^{(k)})\,\mathbf{w}^{(k)}$ and
$\nabla^2 f^{(k)} = -\sin(\phi^{(k)})\,\mathbf{w}^{(k)}(\mathbf{w}^{(k)})^\top$.
Substituting into \eqref{eq:LB_ambient} gives the key formula:
\begin{equation}\label{eq:pw_LB}
  \boxed{
  \Dg f^{(k)}\big|_M
  = -\sin(\phi^{(k)})\,(\mathbf{w}^{(k)})^\top P(x)\,\mathbf{w}^{(k)}
  \;-\; \cos(\phi^{(k)})\,H(x)\cdot\mathbf{w}^{(k)}.}
\end{equation}
The first term involves the \emph{projected squared frequency}
$|\mathbf{w}|_P^2 := \mathbf{w}^\top P(x)\mathbf{w}$, and the second
is a curvature correction proportional to $H(x)\cdot\mathbf{w}$.
Both are computable without automatic differentiation once $P(x)$ and
$H(x)$ are known for the specific manifold.

\begin{remark}[Flat-space recovery]\label{rem:flat}
  For $M = \mathbb{R}^n$ one has $P = I_n$ and $H = 0$, so
  $\Dg f^{(k)} = -|\mathbf{w}^{(k)}|^2\sin\phi^{(k)}$,
  recovering the standard WANPF formula \cite{He2025WANPF}.
\end{remark}

\begin{remark}[Non-eigenfunction character]\label{rem:non_eigen}
  Unlike in flat space, plane waves are \emph{not} eigenfunctions of $\Dg$
  on a curved manifold: the curvature correction introduces a $\cos\phi$ term
  that mixes phases.  Nevertheless, as adversarially optimized test probes,
  they remain effective.
\end{remark}

Using \eqref{eq:pw_LB} and the tangential drift formula
$\langle\mathbf{b}(x), \Ng f^{(k)}(x)\rangle_g
= \cos(\phi^{(k)})\,\langle P(x)\mathbf{b}(x), \mathbf{w}^{(k)}\rangle$,
the backward Kolmogorov operator on a steady-state plane wave evaluates to:
\begin{align}
  \calL_{\mathrm{ss}} f^{(k)}(x)
  &= -\frac{\sigma^2}{2}\sin(\phi^{(k)})\,
  (\mathbf{w}^{(k)})^\top P(x)\mathbf{w}^{(k)} \notag\\
  &\quad - \frac{\sigma^2}{2}\cos(\phi^{(k)})\,
  H(x)\cdot\mathbf{w}^{(k)} \notag\\
  &\quad + \cos(\phi^{(k)})\,\langle P(x)\mathbf{b}(x),\mathbf{w}^{(k)}\rangle.
  \label{eq:Lss_planewave}
\end{align}
All three terms are closed-form at each sample on $M$.

\section{Manifold-Constrained Pushforward Architecture}
\label{sec:arch}

\subsection{Retraction-Based Generator}

The pushforward map must satisfy $F_\vth(t, \mathbf{x}_0, \mathbf{r}) \in M$
for all arguments.  We achieve this by composing an unconstrained residual
network $\tilde{F}_\vth : \mathbb{R}^{1+n+d} \to \mathbb{R}^n$ with a
manifold retraction $\PiM$:
\begin{equation}\label{eq:pf_manifold}
  F_\vth(t, \mathbf{x}_0, \mathbf{r})
  = \PiM\!\Bigl(
    \mathbf{x}_0 + \sqrt{t}\,\tilde{F}_\vth(t, \mathbf{x}_0, \mathbf{r})
  \Bigr).
\end{equation}
Here $\mathbf{x}_0 \sim \rho_0$ is an initial sample on $M$ and
$\mathbf{r} \sim \mathcal{P}_{\mathrm{base}}$ is a $d$-dimensional noise
vector.  The $\sqrt{t}$ prefactor ensures that $F_\vth(0, \mathbf{x}_0,
\mathbf{r}) = \PiM(\mathbf{x}_0) = \mathbf{x}_0$ since $\mathbf{x}_0
\in M$ already, so the initial condition $\rho(0,\cdot) = \rho_0$ is
satisfied automatically.

For commonly encountered manifolds the retraction $\PiM$ has an analytical
form:
\begin{itemize}
  \item \textbf{Sphere} $S^{n-1}$:
    $\PiM(\mathbf{v}) = \mathbf{v}/\|\mathbf{v}\|$.
  \item \textbf{Stiefel manifold} $\mathrm{St}(n,k)$:
    polar retraction $\Pi(\mathbf{V}) = \mathbf{V}(\mathbf{V}^\top\mathbf{V})^{-1/2}$.
  \item \textbf{Flat torus} $\mathbb{T}^n$:
    componentwise wrapping
    $\Pi(\mathbf{v}) = \bigl(\cos v_i, \sin v_i\bigr)_{i=1}^n$.
\end{itemize}

\begin{remark}[Parametrization alternative]
  When $M$ admits a global chart $\psi : \mathbb{R}^m \to M$ (e.g.,
  stereographic projection for $S^{n-1}$, or the exponential map for
  Lie groups), one can instead define
  $F_\vth := \psi(\xi_\vth)$ where $\xi_\vth$ is a network in
  the parameter space $\mathbb{R}^m$.  This avoids the projection step
  but introduces chart singularities for most manifolds of interest;
  the retraction approach \eqref{eq:pf_manifold} is therefore preferred.
\end{remark}

\subsection{Network Specification}

For the sphere experiments in Section~\ref{sec:numerics} we use a
\emph{steady-state} generator (no $t$ or $\mathbf{x}_0$ argument):
\begin{equation}
  F_\vth(\mathbf{r}) = \frac{\tilde{F}_\vth(\mathbf{r})}
  {\|\tilde{F}_\vth(\mathbf{r})\|},
  \qquad
  \tilde{F}_\vth : \mathbb{R}^d \xrightarrow{\text{MLP}} \mathbb{R}^3,
\end{equation}
with $\tanh$ activations throughout.  The network dimensions and training
hyperparameters are listed in Table~\ref{tab:hparams}.

\section{Time-Dependent WANPF on a Manifold}
\label{sec:timedep}

\subsection{Weak Formulation}
\label{sec:timedep_weak}

Let $f \in C^\infty([0,T]\times M)$ be a smooth test function.
Multiplying the time-dependent FPE \eqref{eq:FP_intro} by $f$,
integrating over $M \times [0,T]$, and integrating by parts in both
time and space (boundary terms vanish by compactness) gives:
\begin{align}\label{eq:weak_manifold}
  \int_M f(T,x)\,\rho(T,x)\dVol
  &- \int_M f(0,x)\,\rho_0(x)\dVol \notag\\
  &- \int_0^T\!\!\int_M \rho(t,x)\,
  \Bigl(\partial_t f + \tfrac{\sigma^2}{2}\Dg f + \mathbf{b}(f)\Bigr)\dVol\,dt
  = 0.
\end{align}
Since the pushed-forward samples
$\mathbf{x}(t) = F_\vth(t,\mathbf{x}_0,\mathbf{r})$ are distributed
according to $\rho(t,\cdot)$ when $\vth$ is correct, the identity
\eqref{eq:weak_manifold} becomes three Monte Carlo estimable terms:
\begin{equation}\label{eq:weak_expectation}
  \hat{E}_T^{(k)} - \hat{E}_0^{(k)} - \hat{E}^{(k)} = 0,
\end{equation}
with the time-dependent plane-wave test functions:
\begin{equation}\label{eq:planewave_td}
  f^{(k)}(t,x)
  = \sin\!\bigl(\mathbf{w}^{(k)}\cdot x + \kappa^{(k)} t + b^{(k)}\bigr),
\end{equation}
where $\kappa^{(k)} \in \mathbb{R}$ is a learnable temporal frequency.

\subsection{Monte Carlo Estimators}
\label{sec:mc_estimators}

Let $\phi^{(k)} := \mathbf{w}^{(k)}\cdot x + \kappa^{(k)} t + b^{(k)}$.
The three estimators in \eqref{eq:weak_expectation} are:

\paragraph{Terminal term.}
Draw $M_T$ samples $\mathbf{x}_{0,T}^{(m)}\!\sim\!\rho_0$ and
$\mathbf{r}_T^{(m)}\!\sim\!\mathcal{P}_{\mathrm{base}}$:
\begin{equation}
  \hat{E}_T^{(k)} = \frac{1}{M_T}\sum_{m=1}^{M_T}
  f^{(k)}\!\Bigl(T,\,F_\vth\bigl(T,\mathbf{x}_{0,T}^{(m)},\mathbf{r}_T^{(m)}\bigr)\Bigr).
\end{equation}

\paragraph{Initial term.}
Draw $M_0$ samples $\mathbf{x}_0^{(m)}\!\sim\!\rho_0$ directly:
\begin{equation}
  \hat{E}_0^{(k)} = \frac{1}{M_0}\sum_{m=1}^{M_0}
  f^{(k)}\!\bigl(0,\,\mathbf{x}_0^{(m)}\bigr).
\end{equation}

\paragraph{Interior term.}
Draw $M$ triples
$t^{(m)}\!\sim\!\mathcal{U}(0,T)$,
$\mathbf{x}_0^{(m)}\!\sim\!\rho_0$,
$\mathbf{r}^{(m)}\!\sim\!\mathcal{P}_{\mathrm{base}}$, and set
$\mathbf{x}^{(m)} = F_\vth(t^{(m)},\mathbf{x}_0^{(m)},\mathbf{r}^{(m)})$:
\begin{equation}
  \hat{E}^{(k)} = \frac{T}{M}\sum_{m=1}^{M}
  \Bigl[\kappa^{(k)}\cos(\phi^{(k)})
  + \tfrac{\sigma^2}{2}\Dg f^{(k)}
  + \cos(\phi^{(k)})\,\langle P(x)\mathbf{b}(x),\mathbf{w}^{(k)}\rangle
  \Bigr]_{\!\bigl(t^{(m)},\mathbf{x}^{(m)}\bigr)},
\end{equation}
where $\Dg f^{(k)}$ is evaluated via \eqref{eq:pw_LB}.

\subsection{Adversarial Training Objective}
\label{sec:timedep_loss}

The total loss aggregates the squared residuals over all $K$ test functions:
\begin{equation}\label{eq:loss_td}
  \mathcal{L}_{\mathrm{total}}[\vth,\{\boldsymbol{\eta}^{(k)}\}]
  = \frac{1}{K}\sum_{k=1}^K
  \bigl(\hat{E}_T^{(k)} - \hat{E}_0^{(k)} - \hat{E}^{(k)}\bigr)^2,
\end{equation}
with $\boldsymbol{\eta}^{(k)} = \{\mathbf{w}^{(k)},\kappa^{(k)},b^{(k)}\}$.
Training follows the min-max scheme:
\begin{equation}\label{eq:minmax}
  \min_\vth\;\max_{\{\boldsymbol{\eta}^{(k)}\}}\;
  \mathcal{L}_{\mathrm{total}}.
\end{equation}
The generator (minimiser) and test-function network (maximiser) are updated
with separate Adam optimisers, with gradient ascent for the adversary and
gradient descent (with gradient clipping) for the generator.

\section{Explicit Formulae for Key Manifolds}
\label{sec:manifolds}

We record the projection matrix $P(x)$, the mean-curvature vector $H(x)$,
and the resulting Laplace--Beltrami formula \eqref{eq:pw_LB} for the two
manifolds most frequently encountered in applications.

\subsection{The Unit Sphere \texorpdfstring{$S^{n-1} \subset \mathbb{R}^n$}{Sn-1}}
\label{sec:sphere}

For $x \in S^{n-1}$ (i.e.\ $\|x\| = 1$), the tangential projection and
mean-curvature vector are:
\begin{equation}\label{eq:sphere_geom}
  P(x) = I_n - x x^\top, \qquad H(x) = -(n-1)\,x.
\end{equation}
Substituting into \eqref{eq:pw_LB} gives:
\begin{align}\label{eq:LB_sphere}
  \Dg f^{(k)}\big|_{S^{n-1}}
  &= -\sin(\phi^{(k)})\,
  \bigl(|\mathbf{w}^{(k)}|^2 - (x\cdot\mathbf{w}^{(k)})^2\bigr) \notag\\
  &\qquad + (n-1)\cos(\phi^{(k)})\,(x\cdot\mathbf{w}^{(k)}).
\end{align}
For $S^2 \subset \mathbb{R}^3$ (i.e.\ $n=3$) the constant factor is $n-1 = 2$.

\subsection{The Flat Torus \texorpdfstring{$\mathbb{T}^n \subset \mathbb{R}^{2n}$}{Tn}}
\label{sec:torus}

Embed the $n$-torus as
$M = \{(\cos\theta_i,\sin\theta_i)_{i=1}^n\} \subset \mathbb{R}^{2n}$.
The tangent space is spanned by the $n$ vectors
$e_i = (0,\ldots,-\sin\theta_i,\cos\theta_i,\ldots,0)^\top$, and the
mean-curvature vector vanishes because each factor is a unit circle:
$H(x) = 0$.
Writing $\mathbf{w} = (w_1^c,w_1^s,\ldots,w_n^c,w_n^s)^\top \in \mathbb{R}^{2n}$,
the projected squared frequency is:
\begin{equation}\label{eq:torus_wPw}
  \mathbf{w}^\top P(x)\mathbf{w}
  = \sum_{i=1}^n\bigl(-w_i^c\sin\theta_i + w_i^s\cos\theta_i\bigr)^2,
\end{equation}
and therefore $\Dg f^{(k)} = -\sin(\phi^{(k)})\,\mathbf{w}^\top P(x)\mathbf{w}$,
with no curvature correction.

\section{Numerical Experiment: Double-Well FPE on \texorpdfstring{$S^2$}{S2}}
\label{sec:numerics}

\subsection{Problem Setup}
\label{sec:numerics_setup}

We consider the steady-state Fokker--Planck equation \eqref{eq:FP_steady}
on the unit sphere $S^2 \subset \mathbb{R}^3$ with diffusion coefficient
$\sigma = 0.5$ and potential-derived drift $\mathbf{b}(x) = -P(x)\nabla V(x)$.
The potential is a double-well function on $\mathbb{R}^3$ restricted to $S^2$:
\begin{equation}\label{eq:potential}
  V(x,y,z) = \alpha(x^2 - 1)^2 + \beta z^2,
  \qquad \alpha = 4.0,\; \beta = 2.0.
\end{equation}
The two wells of $V$ are located at $(\pm 1, 0, 0)$ on $S^2$, where $V = 0$.
The $\beta z^2$ term confines mass to the equatorial band and prevents the
distribution from spreading toward the poles.
Figure~\ref{fig:potential} shows the potential $V$ in longitude--latitude
coordinates.

The Gibbs invariant measure is $\rho_\infty(x) \propto e^{-2V(x)/\sigma^2}$,
which concentrates sharply near $(\pm 1,0,0)$ for the chosen parameters.

\subsection{Ambient Gradient and Drift}
\label{sec:numerics_drift}

The Euclidean gradient of $V$ is:
\begin{equation}\label{eq:grad_V}
  \nabla V(x,y,z)
  = \bigl(4\alpha(x^2-1)x,\; 0,\; 2\beta z\bigr)^\top.
\end{equation}
The manifold drift $\mathbf{b}(x) = -P(x)\nabla V(x)$ is then obtained by
projecting $-\nabla V$ onto the tangent plane of $S^2$ at $x$, using
$P(x) = I_3 - xx^\top$.

\subsection{Architecture and Training}
\label{sec:numerics_training}

The generator is a two-layer MLP with width 32 and $\tanh$ activations,
mapping a $d = 3$-dimensional standard Gaussian base distribution
$\mathbf{r} \sim \mathcal{N}(0, I_3)$ to $S^2$ via normalization:
\begin{equation}
  F_\vth(\mathbf{r}) = \frac{\tilde{F}_\vth(\mathbf{r})}
  {\|\tilde{F}_\vth(\mathbf{r})\|}.
\end{equation}
The adversarial test functions are $K = 200$ ambient plane waves
\eqref{eq:planewave_ss} with parameters initialized as
$\mathbf{w}^{(k)} \sim \mathcal{N}(0, 4I_3)$ and $b^{(k)} \sim \mathcal{N}(0, 0.25)$.
The $\calL_{\mathrm{ss}} f^{(k)}(x)$ evaluator uses the closed-form
formula \eqref{eq:Lss_planewave} specialized to $S^2$ via \eqref{eq:LB_sphere},
with no automatic differentiation through the geometry.
Training hyperparameters are summarized in Table~\ref{tab:hparams}.

\begin{table}[ht]
  \centering
  \caption{Hyperparameters for the $S^2$ double-well experiment.}
  \label{tab:hparams}
  \begin{tabular}{lll}
    \toprule
    Parameter & Symbol & Value \\
    \midrule
    Diffusion coefficient    & $\sigma$       & $0.5$ \\
    Base dimension           & $d$            & $3$ \\
    Generator hidden width   & ---            & $32$ \\
    Generator depth          & ---            & $2$ layers \\
    Number of test functions & $K$            & $200$ \\
    Monte Carlo batch size   & $M_s$          & $200$ \\
    Training steps           & $N$            & $5{,}000$ \\
    Generator learning rate  & $\eta_{\vth}$  & $10^{-3}$ (cosine annealing) \\
    Adversary learning rate  & $\eta_{\eta}$  & $5\times10^{-3}$ \\
    Potential parameters     & $\alpha,\beta$ & $4.0,\;2.0$ \\
    \bottomrule
  \end{tabular}
\end{table}

The training loop alternates one adversary gradient-ascent step (maximizing
the loss over $\{\boldsymbol{\eta}^{(k)}\}$) with one generator
gradient-descent step (minimizing the loss over $\vth$), with gradient
clipping of norm $1.0$ applied to the generator.
A cosine annealing schedule is applied to the generator learning rate
from $10^{-3}$ down to $10^{-5}$ over the $5{,}000$ steps.

\subsection{Results}
\label{sec:numerics_results}

\begin{figure}[htb]
  \centering
  \includegraphics[width=\textwidth]{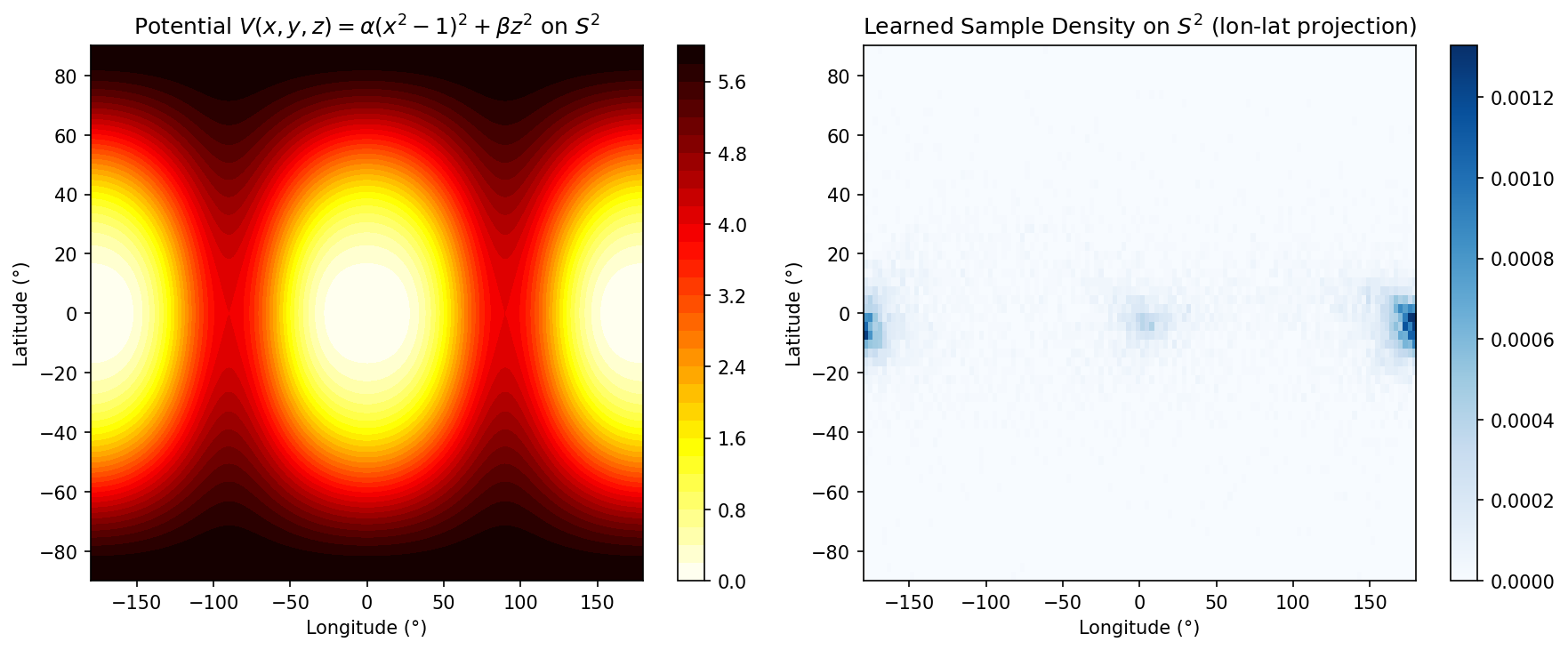}
  \caption{Left: the double-well potential $V(x,y,z) = \alpha(x^2-1)^2 + \beta z^2$
    on $S^2$, shown in longitude--latitude projection.
    The two potential minima ($V = 0$) are located at longitude $0^\circ$ and
    $\pm 180^\circ$ on the equator, corresponding to $(\pm 1, 0, 0)$.
    Right: density histogram of $8{,}000$ samples generated by the trained
    neural pushforward map in the same coordinate system.
    Mass concentrates near both well bottoms and is absent near the poles,
    consistent with the Gibbs distribution
    $\rho_\infty \propto e^{-2V/\sigma^2}$ for $\sigma = 0.5$.}
  \label{fig:potential}
\end{figure}

Figure~\ref{fig:potential} shows the double-well potential $V$ on $S^2$ in
a longitude--latitude (Plate Carrée) projection alongside the histogram
of learned sample density produced by the trained generator.

\begin{figure}[htb]
  \centering
  \includegraphics[width=\textwidth]{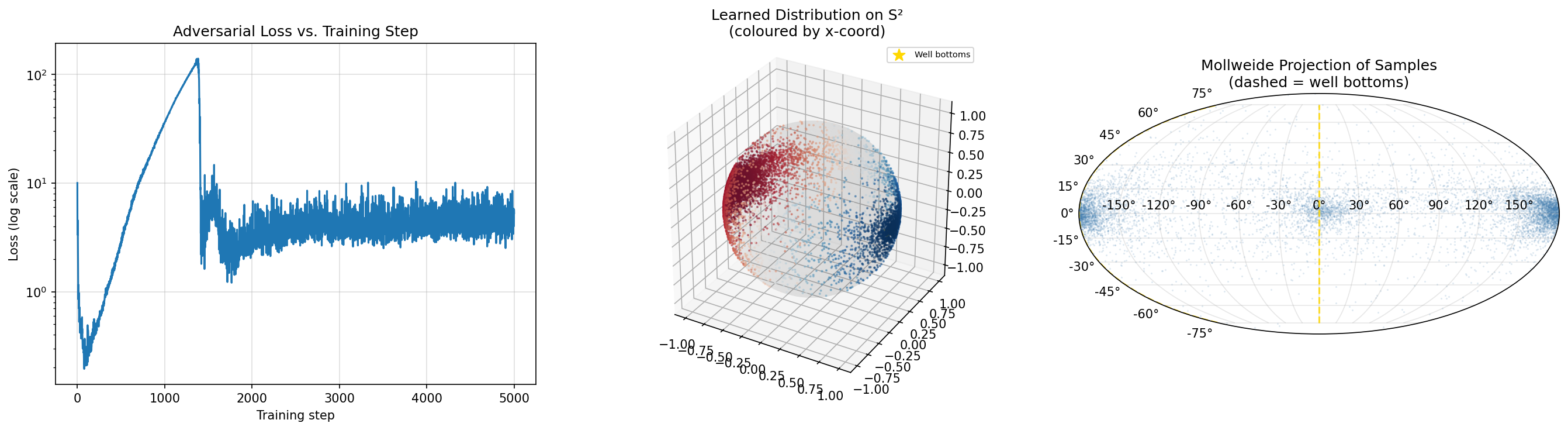}
  \caption{Training and sampling diagnostics for the $S^2$ double-well experiment.
    Left: adversarial loss versus training step (log scale).
    The loss rises during early generator exploration and then decreases as
    training converges to a stationary min-max point near step $2{,}000$.
    Center: $8{,}000$ generated samples on $S^2$, coloured by their
    $x$-coordinate value (red near $x=+1$, blue near $x=-1$); gold stars
    mark the two well bottoms $(\pm 1, 0, 0)$.
    The generator correctly concentrates mass around both minima symmetrically.
    Right: Mollweide equal-area projection of the same samples; dashed gold
    lines indicate the longitudes of the well bottoms.
    Samples cluster near $0^\circ$ and $\pm 180^\circ$ longitude, confirming
    two-well capture, with near-uniform spreading in longitude about each minimum
    and strong equatorial confinement consistent with the $\beta z^2$ term.}
  \label{fig:results}
\end{figure}

Figure~\ref{fig:results} shows additional diagnostics: the adversarial
training loss on a log scale, the pushed-forward samples on the sphere
coloured by the $x$-coordinate, and a Mollweide projection of the samples.

The key qualitative observations are as follows.
The learned distribution concentrates mass near both well bottoms at $(\pm 1, 0, 0)$
and is strongly suppressed near the poles, as expected from the $\beta z^2$
confinement term in $V$.  The approximate symmetry between the two wells is
preserved, reflecting the $x \to -x$ symmetry of $V$.  The training loss
converges to a stationary level after approximately $2{,}000$ steps, with
subsequent fluctuations characteristic of the min-max adversarial dynamics.

\section{Discussion}
\label{sec:discussion}

\subsection{Advantages of the Manifold Formulation}

The manifold WANPF inherits all structural advantages of the flat-space
method \cite{He2025WANPF}: samples are generated rather than a pointwise
density, so probability conservation is exact; the weak formulation requires
no mesh, no chart atlas, and no Jacobian determinant; and the method handles
multimodal distributions without any explicit mixture structure.

The manifold extension adds two further benefits.  First, the retraction
architecture $F_\vth = \Pi_M \circ \tilde{F}_\vth$ keeps all samples on $M$
at no additional training cost beyond a single projection per forward pass.
Second, the closed-form Laplace--Beltrami formula \eqref{eq:pw_LB} bypasses
all automatic differentiation through the geometry of $M$, making the
per-step cost independent of the number of layers in the intrinsic
differential operators.

\subsection{Connection to Geometric Measure Theory}

The use of ambient test functions is consistent with the theory of varifolds,
where the key identity
\begin{equation}
  \int_M \rho\,\Dg f\dVol
  = \int_M \rho\,\bigl(\nabla^2 f : P - \nabla f\cdot H\bigr)\dVol
\end{equation}
expresses the Laplace--Beltrami action entirely through ambient quantities
evaluated on samples in $\mathbb{R}^n$ \cite{simon1983lectures}.
WANPF exploits precisely this identity.

\subsection{Extension to Non-Isotropic Diffusion}

The derivation extends directly to the general FPE with a diffusion tensor
$a$ (a positive-definite $(2,0)$-tensor on $M$).  The backward Kolmogorov
operator becomes
$\calL f = \tfrac{1}{2}\Eg(a(\mathrm{d}f)) + \mathbf{b}(f)$,
and in ambient coordinates:
\begin{equation}\label{eq:anisotropic_L}
  \calL f = \frac{1}{2}\,\nabla^2 f : P A P^\top
  - \frac{1}{2}(A P^\top \nabla f)\cdot H
  + \langle P\mathbf{b},\, P\nabla f\rangle,
\end{equation}
where $A \in \mathbb{R}^{n\times n}$ is any ambient extension of $a$.
For the isotropic case $a = \sigma^2 g^{-1}$, $P A P^\top = \sigma^2 P$
and \eqref{eq:anisotropic_L} reduces to \eqref{eq:Lss_planewave}.

\subsection{Supplementing with Spherical Harmonics}

As noted in Remark~\ref{rem:non_eigen}, ambient plane waves are not
eigenfunctions of $\Dg$ on a curved manifold.  For $S^{n-1}$, one could
supplement the plane-wave test-function class with spherical harmonics
$Y_\ell^m$, which \emph{are} eigenfunctions with eigenvalue
$-\ell(\ell + n - 2)$.  In the WANPF framework this amounts to adding
additional test functions whose operator evaluations are analytically known;
the min-max training structure is unchanged.  We leave a systematic
comparison to future work.

\subsection{Future Directions}

Several extensions are natural.  The time-dependent formulation
(Section~\ref{sec:timedep}) is directly applicable to the heat equation on
$S^{n-1}$ or on a Lie group such as $\mathrm{SO}(3)$; the only change is
the specification of $P(x)$ and $H(x)$ for the new manifold.
For the fractional Fokker--Planck equation on a manifold, one would need to
replace the Laplace--Beltrami operator with the spectral fractional
Laplacian $(-\Dg)^{\alpha/2}$, whose plane-wave action on curved manifolds
is no longer as clean but could be approximated via Monte Carlo integration
over manifold heat kernels.  Higher-dimensional manifolds (e.g., $S^{n-1}$
with $n \geq 10$) are a natural stress test for the mesh-free approach.

\section{Conclusion}
\label{sec:conclusion}

We have extended the Weak Adversarial Neural Pushforward Method to the
Fokker--Planck equation on compact embedded Riemannian manifolds.
The extension rests on three pillars: (i) a retraction-based pushforward
architecture that keeps samples on $M$ by construction; (ii) ambient
plane-wave test functions whose Laplace--Beltrami operators are computable
in closed form via the tangential projection $P(x)$ and mean-curvature
vector $H(x)$; and (iii) a min-max adversarial loss whose Monte Carlo
estimator requires only pointwise evaluations at samples on $M$, with no
mesh, no chart, and no Jacobian computation.

The resulting algorithm for the steady-state case is summarized in
Algorithm~\ref{alg:wanpf_steady}.
Numerical experiments on the double-well FPE on $S^2$ confirm that the
trained generator correctly concentrates mass near both potential minima,
consistent with the Gibbs invariant measure.

\begin{algorithm}[htb]
\caption{Steady-State Manifold WANPF}
\label{alg:wanpf_steady}
\begin{algorithmic}[1]
  \Require manifold $M$ with retraction $\Pi_M$, $P(\cdot)$, $H(\cdot)$;
    drift $\mathbf{b}$; diffusion $\sigma$;
    learning rates $\eta_\vth$, $\eta_\eta$; batch size $M_s$; steps $N$.
  \State Initialize generator $F_\vth$ and test-function parameters
    $\{\boldsymbol{\eta}^{(k)}\}_{k=1}^K$.
  \For{step $= 1, \ldots, N$}
    \State \textbf{Adversary step:} sample $\{\mathbf{r}^{(m)}\}$, compute
      $\mathbf{x}^{(m)} = F_\vth(\mathbf{r}^{(m)}) \in M$, evaluate
      $\hat{E}^{(k)}$ via \eqref{eq:Ehat_steady} and \eqref{eq:Lss_planewave},
      update $\boldsymbol{\eta}^{(k)} \mathrel{+}= \eta_\eta \nabla_{\boldsymbol{\eta}}
      \mathcal{L}$ (gradient ascent).
    \State \textbf{Generator step:} same forward pass, update
      $\vth \mathrel{-}= \eta_\vth \nabla_\vth \mathcal{L}$ (gradient descent).
  \EndFor
  \Ensure generator $F_\vth$ approximating the invariant measure $\rho_\infty$.
\end{algorithmic}
\end{algorithm}

The framework is directly applicable to any manifold for which $P(x)$ and
$H(x)$ can be computed analytically or numerically, and it scales to
dimensions inaccessible by mesh-based solvers.

\bibliographystyle{plain}

\end{document}